\title{Positional games on random graphs}
\author{Milo\v s Stojakovi\'c
\thanks{Institute of Theoretical Computer Science, ETH Zurich, CH-8092
  Switzerland. 
Email addresses: $\{$smilos, szabo$\}$@inf.ethz.ch} \thanks{Supported by the
  joint Berlin/Zurich graduate
  program Combinatorics, Geometry and Computation,
  financed by ETH Zurich and the German Science
  Foundation (DFG).}
\and
Tibor Szab\'o\setcounter{footnote}{0}\footnotemark
}
\documentclass[11pt]{article} 

\usepackage{a4}
\usepackage{latexsym}
\usepackage{amssymb}

\newtheorem{theorem}{Theorem}
\newtheorem{lemma}[theorem]{Lemma}
\newtheorem{proposition}[theorem]{Proposition}
\newtheorem{corollary}[theorem]{Corollary}
\newtheorem{definition}{Definition}
\newtheorem{conjecture}{Conjecture}
\newtheorem{problem}{Problem}

\newcommand\sep{:\,}

\newcommand\ExpOp{\hbox{\bf E}}
\newcommand{\Ex}[1]{\ExpOp\hspace{-0.3ex}\left[#1\right]}
\newcommand{\Prob}[1]{{\rm Pr}\hspace{-0.2ex}\left[ #1\right]}

\def\proof{\noindent{\bf Proof. \,\, }}

\def\cupdis{\dot \cup}
\def\eps{\varepsilon}
\def\:{\colon}

\def\cf{{\cal F}}

\def\cm{{\cal M}}
\def\ct{{\cal T}}
\def\ch{{\cal H}}
\def\cp{{\cal P}}
\def\ckk{{\cal K}}

\def\eopf{\proofend}

\makeatletter
\newcommand{\ProofEndBox}{{\ifhmode\unskip\nobreak\hfil\penalty50 \else
          \leavevmode\fi\quad\vadjust{}\nobreak\hfill$\Box$
            \finalhyphendemerits=0 \par}}%
\makeatother

\newcommand{\proofend}{\ProofEndBox\smallskip}

\newcommand{\heading}[1]{\vspace{1ex}\par\noindent{\bf #1}}

\begin{document}
\maketitle

\begin{abstract}
We introduce and study Maker/Breaker-type positional games on
random graphs. 
Our main concern is to determine the threshold
probability $p_{\cf}$ for the existence of Maker's strategy 
to claim a member of $\cf$ in the unbiased game played on the edges of
random graph $G(n,p)$, 
for various target families $\cf$ of winning sets.
More generally, for each probability above this threshold
we study the smallest
bias $b$ such that Maker wins the $(1\:b)$ biased game.
We investigate these functions for a number of basic games,
like the connectivity game, the perfect matching game, the clique game
and the Hamiltonian cycle game.
\end{abstract}

\section{Introduction} \label{s:intro}

\paragraph{(Un)biased positional games.}
Let $X$ be a finite nonempty set and $\cf\subseteq 2^X$. The pair
$(X, \cf)$ is a {\em positional game} on $X$.
The game is played by two players Maker and Breaker, where
in each move Maker claims one previously unclaimed element of
$X$ and then Breaker claims one previously unclaimed element of
$X$. Maker wins if he claims all
the elements of some set in $\cf$, otherwise Breaker wins.
The set $X$ will be referred to as the board, 
and the set $\cf$ as the set of winning sets. Whenever there is
no confusion about what the board is, 
we may refer to the game $(X,\cf)$ as just $\cf$.

Unless otherwise stated, we assume that Maker starts the game.
We note, however, that the asymptotic statements discussed in the paper
are not influenced by which player makes the first move. 
For technical reasons we still have to talk about games
in which Breaker starts. So in order to avoid confusion, the positional
game
with board $X$ and set of winning sets $\cf$
in which Breaker makes the first move is 
denoted by $(\widehat X, \cf)$.

The set of all positional games could be partitioned into two
classes.
The game $(X,\cf)$ is called a {\em Maker's win} if 
Maker has a winning strategy, that is, 
playing against an arbitrary strategy Maker can 
occupy a member of $\cf$. Clearly, if $(X,\cf )$ is {\em not} a 
Maker's win, then Breaker is able to
prevent any opponent from occupying a 
winning set. Such a positional game is called a {\em Breaker's win}.

Typical, well-studied examples of such positional games are played
on the edges of a complete graph, i.e.\ $X=E(K_n)$. 
Maker's goal usually is to build a graph theoretic structure -- like a 
spanning tree, a perfect matching, a Hamiltonian cycle, or a clique of 
fixed size.
It turns out that all these games are won easily by Maker if $n$ is 
sufficiently big, so in order
to make things more fair (if such thing exists; actually no game of
perfect information is {\em fair} as the winner---in theory---is known in
the beginning of the game) one could give Breaker extra power by allowing
him to claim more than $1$ edge in each move.

If $X$ is a finite nonempty set, $\cf\subseteq 2^X$ and $a,b$ are positive
integers (possibly functions of the board size), then the 4-tuple
$(X, \cf, a, b)$ is a {\em biased $(a\:b)$ game}.
In a biased $(a\:b)$ game, Maker claims $a$ elements (instead of 1) 
and Breaker claims $b$ elements (instead of 1) in each move.
Recall, that unless otherwise stated Maker starts the game. 
The biased game in which Breaker starts
is denoted by $(\widehat X, \cf, a, b)$. Note that $a$ is always 
the bias of Maker, independently from who is the first player to move. 

For a family $\cf$ the smallest integer $b_{\cf}$ is sought 
(and sometimes found; see \cite{bec, be2, be3, bl1, bl2, ce}) 
for which Breaker wins the $(1:b_{\cf})$ game.

In the {\em connectivity game} Maker's goal is to build a connected 
spanning subgraph; i.e. in this game the family of winning sets
is the family $\ct=\ct_n$ of all spanning trees on $n$ vertices. 
Chv\' atal and Erd\H os proved \cite{ce}
that $b_\ct=\Theta (\frac{n}{\log n} )$.

Beck \cite{bec} established $b_\ch=\Theta(\frac{n}{\log n})$, where
$\ch=\ch_n$ is the family of all Hamiltonian cycles on $n$ vertices.

For the family $\ckk_k=\ckk_{k,n}$
of all $k$-cliques on $n$ vertices,  
Bednarska and \L uczak \cite{bl1} 
showed that
$b_{\ckk_k}=\Theta (n^{\frac{2}{k+1}})$. More generally, they proved
that in the game in which Maker's goal is to claim an arbitrary fixed
graph $G$, the threshold bias is $\Theta (n^{1/m'(G)})$. (Here $m'(G)$ is
the maximum of $\frac{e(H)-1}{v(H)-2}$ over all subgraphs $H$ of $G$ with
at least $3$ vertices.)

\paragraph{Playing on a random board.}
In the present paper we introduce another approach to even out the
advantage Maker has in a $(1\: 1)$ game, by randomly reducing the board
size and keeping only those winning sets which survive this
thinning intact. 

\begin{definition} \label{df:1}
Let $(X, \cf, a, b)$ be a biased game. {\em Random game}
$(X_p, \cf_p,a,b)$ is a probability space of games  where each $x\in X$
is independently included in $X_p$ with probability $p$, 
and $\cf_p=\{W\in \cf \sep W\subseteq X_p\}$.
\end{definition}

Apart from the trivial case $\emptyset \in \cf$,
Breaker surely wins when $p=0$.
On the other hand, the unbiased version of all the 
graph games that we consider are (easy) Maker's wins,
when $p=1$ and the board is sufficiently large.
For any other probability $p$, $0<p<1$,
we cannot be sure who (Maker or Breaker) wins the random 
game $\cf_p$. The best we can
conclude is that Maker (or Breaker) wins a.s.\ (almost
surely), i.e.\ the probability that Maker (Breaker) wins tends
to 1 if the board size tends to infinity. (So we actually talk about an
infinite family of probability spaces of games \dots)

Let $(X,\cf)$ be a particular sequence of games, where 
$\emptyset\notin \cf$, the board size tends to infinity, 
and $(X,\cf,1,1)$ is 
won by Maker provided $|X|$ is big enough. The first natural
question to ask is: What is the threshold probability
$p_{\cf}$ at which an almost sure Breaker's win turns into an
almost sure Maker's win. More precisely we would like to determine 
$p_\cf$ for which 
\begin{itemize}
\item $\Prob{(X_p,\cf_p,1,1) \mbox{ is a Breaker's win}}\rightarrow 1$ for 
$p=o(p_{\cf})$, and \\
\item $\Prob{(X_p,\cf_p,1,1) \mbox{ is a Maker's win}}\rightarrow 1$ 
for $p=\omega(p_{\cf})$.
\end{itemize}

Such a threshold $p_\cf$ exists \cite{bt}, since being a Maker's win 
is an {\em increasing property}.

The main goal of this paper is to establish a connection between the 
natural threshold values, $b_\cf$ and $p_\cf$, 
corresponding to the two different weakenings of Maker's power:
bias and random thinning, respectively.
We find that there is an intriguing reciprocal connection between these two
thresholds in a number of well-studied games on graphs.

Recall the notations $\ct$, $\ch$, and $\ckk_k$, and let us 
denote by $\cm$ the set of all perfect matchings on the graph $K_n$.

\begin{theorem}\label{main1}
For positional games, played on $E(K_n)$, we have
\begin{itemize}
\item[(i)] $p_\ct = \frac{\log n}{n} $, \\
\item[(ii)] $p_\cm = \frac{\log n}{n}$, \\
\item[(iii)] $\frac{\log n}{n} \leq p_\ch \leq 
\frac{\log n}{\sqrt{n}}$, \\
\item[(iv)] $n^{-\frac{2}{k+1}-\varepsilon}  \leq 
p_{\ckk_k} \leq  n^{-\frac{2}{k+1}}$,
for every integer $k\geq 4$ and every constant $\varepsilon>0$.
\\
\item[(v)] 
$p_{\ckk_3} =  n^{-\frac{5}{9}}$.
\end{itemize}
\end{theorem}

For the connectivity game $\ct$ an even more precise statement is true.
In Corollary~\ref{connectivity-precise} we observe that Maker starts
to win a.s.\ at the very moment when the last vertex of a random graph
process picks up its second incident edge.

More generally, for every $p$ we would like to 
find the smallest bias
$b^p_{\cf}$ such that 
Breaker wins the random game $(X_p, \cf_p, 1, b^p_\cf)$ 
a.s. Note that by definition $b_\cf=b^1_\cf$. Another trivial
observation is that $b^p_\cf=0$ 
provided $p$ is less than the threshold for the appearance of
the first element of $\cf$ in the random graph.

We obtain the following.
\begin{theorem}\label{main2}
There exist constants $C_1, C_2, C_3$, such that
\begin{itemize}
\item[(i)] $b_\ct^p  = \Theta \left(pb_\ct\right) = 
\Theta \left( p\frac{n}{\log n}\right)$, provided $p\geq
C_1\frac{1}{b_{\ct}}$,\\
\item[(ii)] $b_\cm^p = \Theta \left(p b_\cm\right) = 
\Theta \left(p\frac{n}{\log n}\right)$, provided $p\geq 
C_2\frac{1}{b_{\cm}}$,\\
\item[(iii)] $\Omega\left(p\frac{\sqrt{n}}{\log n}\right)
\leq b_\ch^p \leq O \left(p\frac{n}{\log n}\right) $, provided $p\geq
C_3\frac{\log n}{\sqrt{n}}$, \\
\item[(iv)] There exists $c_k>0$, such that 
$b_{\ckk_k}^p = \Theta\left( pb_{\ckk_k}\right) 
= \Theta \left( p n^{\frac{2}{k+1}} \right)$,  
provided $p= \Omega\left(\frac{\log^{c_k} n}{b_{\ckk_k}}\right)$.
\end{itemize}
\end{theorem}

One can see that $b_\cf^p$ is of order $p/p_\cf=pb_{\cf}$ 
for the connectivity game and the perfect 
matching game, provided $p\geq Cp_\cf$ for some constant $C$. 
In particular for these games $p_\cf=\Theta(1/b_\cf )$.
In part $(iv)$ of Theorem~\ref{main2}, 
generalizing the arguments of Bednarska and \L uczak
\cite{bl1} we show that one can estimate $b_{\ckk_k}^p$ up to a
constant factor, for all probabilities down to 
a polylogarithmic factor away from the critical probability 
$1/b_{\ckk_k}=n^{-\frac{2}{k+1}}$.
On the other hand Theorem~\ref{main1} part $(v)$ 
shows that in the case $k=3$ we cannot get arbitrarily close to probability
$1/b_{\ckk_k}$,
since Maker
{\em can win} even for probabilities below $1/b_{\ckk_3}=n^{-1/2}$.

Nevertheless we think the Hamiltonian cycle game behaves ``nicely'', 
i.e. the same way
as the connectivity game and the perfect matching game.
\begin{conjecture}
Let $\ch$ be the set of Hamiltonian cycles in $K_n$. There exists a 
constant $C$ such that
$$
b_{\ch}^p = \Theta \left( p\frac{n}{\log n}\right), \mbox{ provided
$p\geq C\frac{\log n}{n}$}.$$
In particular,
$$
p_{\ch} = \frac{\log n}{n}.
$$
\end{conjecture}

Observe that the validity of the conjecture would mean that in a 
random graph with edge probability
$p\geq C\frac{\log n}{n}$ 
Maker could build a Hamiltonian cycle. 
So P\'osa's Theorem (which only proves the existence
of a Hamiltonian cycle) would be true constructively even if an
adversary is playing against us.

The paper is organized as follows. In
Section \ref{s:criterion} we prove a general criterion for Breaker's win
in a different, auxiliary random game. 
In Section~\ref{s:games}, the analysis of
four biased random games is presented.
In particular, in Subsections~\ref{ss:conn},
\ref{ss:ham}, \ref{ss:pm} and \ref{ss:clique} we
look at the connectivity game, the Hamiltonian cycle game, the perfect
matching game and the clique game, respectively.
In Section~\ref{s:unbiased} we analyze more precisely a couple
of $(1\:1)$ games --
the connectivity game (Subsection~\ref{ss:con1on1}) and the clique
game (Subsection~\ref{ss:clique1on1}). Finally, in 
Section~\ref{s:open} we
give a collection of open questions and conjectures.

\heading{Notation.} For a graph $G$, $e(G)$ and $v(G)$ denote the
number of edges and vertices (respectively) of $G$, $\delta (G)$
denotes the minimum degree of $G$, and $E(G)$ and
$V(G)$ denote the sets of edges and vertices (respectively). If
$C\subseteq V(G)$ and $v\in V(G)$, then 
$N_C(v)$ denotes the set of neighbors of $v$ in $C$.
The logarithm $\log n$ in this paper is always of natural base.
For functions $f(n), g(n)\geq 0$, we say that $f=O(g)$ if there are
constants $C$ and $K$, such that $f(n) \leq Cg(n)$ for $n\geq K$;
$f=\Omega( g)$ if $g=O(f)$; 
$f=\Theta (g)$ if $f=O(g)$ and $f=\Omega (g)$;
$f=o(g)$ if $f(n)/g(n)\rightarrow 0$ when $n\rightarrow \infty$;
$f=\omega (g)$ if $g=o(f)$.

\section{A criterion} \label{s:criterion}

One of few general, but still very applicable results 
to decide the winner of biased positional games is the biased 
version of the Erd\H os--Selfridge Theorem \cite{es, be2}. 
It provides a criterion for Breaker to win, 
applicable on any game.
\begin{theorem} (Beck, \cite{be2}) \label{t:es}
If 
$$
\sum_{A\in \cf} (1+b)^{-|A|/a} < 1,
$$ 
then Breaker has a winning strategy in the $(\widehat{X},\cf,a,b)$ game.
\end{theorem}
If Maker plays the first move then the $1$ 
on the right hand side of the criterion is to be replaced by 
the fraction $\frac{1}{1+b}$.

We will also need the following extension.

\begin{theorem} (\cite{be2, bl1}) \label{t:gen-es}
If for a positive integer $c$ we have
$$
\sum_{A\in \cf} (1+b)^{-|A|/a} < c\frac{1}{1+b},
$$ 
then Breaker has a winning strategy in the 
$(X,\{\cup_{B\in F} B \sep F\in {\cf \choose c} \},a,b)$ game.
\end{theorem}

In this section we give an adaptation of the first 
criterion which proves to be 
very useful in dealing with positional games on a random board.
We need the following technical definition.
\begin{definition}
Let $(X, \cf, a, b)$ be a biased game. 
{\em Random game $(X_p, \cf_p^\cap,a,b)$ with induced set of winning sets} 
is a probability space of games, where $X_p$ 
is defined as in Definition~\ref{df:1} and
$\cf_p^\cap=\{W\sep \exists F\in \cf,\, W=F\cap X_p\}$.
\end{definition}

The following statement is the randomized version of Theorem~\ref{t:es}.
It is stated for the biased $(b\:1)$ game in which Breaker
is the first player, because this is the version we will need
in our applications. 

\begin{theorem} \label{t:es_cap}
Let $\cf$ be a set of winning sets on $X$ with
\begin{eqnarray}
\sum_{A\in\cf} 2^{-\frac{|A|}{b}} <1 \label{c:es}
\end{eqnarray}
(i.e. the condition of the Erd\H os--Selfridge Theorem holds for the
$(\widehat X, \cf, b, 1)$ game), and
\begin{eqnarray}
\lim_{n\rightarrow\infty} \min_{A\in\cf} \frac{|A|}{b} =
\infty. \label{c:ws_large}
\end{eqnarray}

If $p$ and $\delta >0$ are chosen so that
$p>\frac{4\log 2}{\delta^2 b}$ holds, then
the game $(\widehat{X}_p, \cf_{p}^\cap,
(1-\delta)pb, 1)$ is a Breaker's win a.s.
\end{theorem}

\proof
For each $A\in\cf$ and its corresponding set
${A'} \in\cf_{p}^\cap$ we have $\Ex{|{A'}|} = p|A|$.
If all winning sets $A'\in\cf_{p}^\cap$ have size at least
$(1-\delta) p|A|$, then
$$
\sum_{{A'}\in \cf_{p}^\cap} 2^{-\frac{|A'|}{(1-\delta)pb}}
\leq
\sum_{{A}\in \cf} 2^{-\frac{(1-\delta)p|A|}{(1-
\delta)pb}} =
\sum_{{A}\in \cf} 2^{-\frac{|A|}{b}}
<1.
$$
Using the Erd\H os--Selfridge theorem we obtain that Breaker
wins the $(\widehat{X}_p, \cf_{p}^\cap,(1-\delta)pb, 1)$ 
game, 
provided $|A'|\geq (1-\delta) p|A|$ for all $A'\in \cf_{p}^\cap$.

Next we check that this condition holds almost surely.
Using a Chernoff bound, we obtain that
$$
\Prob{\exists A\in\cf \sep |A'|\leq (1-\delta) p|A|} \leq
\sum_{A\in\cf}
e^{-\frac{\delta^2p |A|}{2}}.$$

If we denote $\min_{A\in\cf} \frac{|A|}{b} $ by $m_n$, then we have
\begin{eqnarray*}
\sum_{A\in\cf}
e^{-\frac{\delta^2p |A|}{2}} \leq
\sum_{A\in\cf}
2^{-2\frac{|A|}{b}} \leq \sum_{A\in\cf}
2^{-m_n} 2^{-\frac{|A|}{b}} <
2^{-m_n} \rightarrow 0,
\end{eqnarray*}
and therefore all winning sets $A'\in\cf_{p}^\cap$ have size at least
$(1-\delta) p|A|$ a.s.

\eopf

\section{Games} \label{s:games}

\subsection{Connectivity game} \label{ss:conn}

The first game we study is a random version of the  
biased connectivity game $(E(K_n), \ct, 1,b)$ on 
a complete graph on $n$ vertices $K_n$. Maker's goal is to build
a spanning, connected subgraph, i.e.\ $\ct$ is the set of
all spanning trees on $n$ vertices.

It is obvious that $p_{\ct}=\Omega (\frac{\log n}{n})$, since
for lower probabilities the random graph is a.s.\ not connected, and Breaker
wins even if he does not claim any edges.

First we generalize this for arbitrary probability $p$ by providing
Breaker with a strategy to isolate a vertex.
One of our main tools is the following winning criterion of
Chv\' atal and Erd\H os on games with disjoint winning sets.

 \begin{theorem} \label{t:C-E} \cite{ce}
 In a biased $(b\:1)$
 game with $k$ disjoint winning sets of size $s$ Maker wins if
\begin{equation}\label{e:C-E}
 s\leq (b-1)\sum_{i=1}^{k-1} \frac{1}{i}.
\end{equation}
 \end{theorem}

 \begin{corollary} \label{c:C-E}
 In a biased $(b\:2)$
 game with $k+1$ disjoint winning sets of size at most $s$ Maker wins if
 $$
 s\leq \left(\left\lfloor\frac{b}{2}\right\rfloor-1\right)\sum_{i=1}^{k-1} 
 \frac{1}{i}.
 $$
 \end{corollary}

\noindent
{\bf Proof of Corollary.} Recall that as a default Maker starts the 
game in the Theorem and the Corollary as well. 
Now Theorem~\ref{t:C-E} obviously remains true
(i.e. Maker wins) even if Breaker starts, provided there are $k+1$
disjoint winning sets instead of $k$. This implies that when 
Breaker starts, the bias is $(2b:2)$, there are $k+1$ winning sets 
and (\ref{e:C-E}) holds, then Maker still wins. Indeed,
since the winning sets are disjoint, after Breaker's
move Maker can 
just pretend to play a $(b:1)$ game and answer with his first $b$ moves
to one of the two selections of Breaker, and answer with his second $b$
moves to
the other move of Breaker, both according to the $(b:1)$ strategy.
Now the Corollary follows, since starting instead of being second
player cannot hurt Maker.
\eopf

 \begin{theorem} \label{t:cg-b}
 There exists $K_0>0$ so that for arbitrary $p\in [0,1]$ and
 $b\geq K_0 p \frac{n}{\log n}$
 Breaker, playing the $(1:b)$ game on the edges of random graph $G(n,p)$,
 can achieve that Maker's graph has an isolated vertex a.s.
 \end{theorem}

 \proof
 Let us fix $b=\lfloor K_0pn/\log n \rfloor $, 
 where $K_0$ is a constant to be determined later. Note that we can assume 
 $p>\log n/2n$, since otherwise the random graph does have an isolated vertex
 a.s., thus Breaker achieves his goal without having to play any 
 moves.

 We present a strategy for Breaker to claim all the edges incident to
 some vertex of $G(n,p)$. If successful, this strategy prevents
 Maker from building a connected subgraph. Similar strategy was
 introduced by Chv\' atal and Erd\H os \cite{ce} for solving the
 problem on the complete graph.

Let $C$ be an arbitrary subset of the
vertex set of cardinality $\lfloor 
n/\log n \rfloor$. Breaker will claim all the edges
incident to some vertex $v\in C$ (thus preventing Maker from claiming 
any edge incident to $v$). 
We would like to use the game from Corollary~\ref{c:C-E}, with the
winning sets being the $\lfloor n/\log n\rfloor$
stars of size at most $n-1$ whose center is in $C$.
Since these stars are not necessarily 
disjoint, formally we will talk about 
ordered pairs of vertices: the winning sets are denoted by 
$W_v=\{ (v, u): u\in V\}$, $v\in C$. We call this game {\em Box}.
To avoid confusion with Maker and Breaker of the game from
Theorem~\ref{t:cg-b}, the players from Corollary~\ref{c:C-E} will 
be called
BoxMaker and BoxBreaker. Recall that in Box the bias is $(b:2)$.

Breaker will utilize the strategy of BoxMaker from Corollary~\ref{c:C-E}
to achieve his goal. How? He will play a game of Box in such a way that
a win for BoxMaker automatically implies a win for Breaker.
When Maker selects an edge $uv$, Breaker interprets it as BoxBreaker
claimed the elements $(u,v)$ and $(v,u)$ in Box. 
Whenever Breaker would like to make a 
move, he looks at the current move of BoxMaker in Box, and takes
those edges which correspond to the $b$ ordered pairs BoxMaker
selected.
If he is supposed to select an edge
which has already been 
selected by him, he selects an arbitrary unoccupied edge. Note that the
above strategy never calls for Breaker to select an edge which 
has already been selected by Maker.

It is also obvious, that if BoxMaker wins Box, then Breaker occupied 
all incident edges of a vertex from $C$.

In order to apply Corollary~\ref{c:C-E} it is enough then to show that the 
 size $d(v)$ of each winning set is appropriately bounded from above, i.e.
 for each $v\in C$ we have 
 $d(v) \leq \frac{K_0}{8} pn\leq 
 \left(\left\lfloor\frac{b}{2}\right\rfloor-1\right) \sum_{i=1}^{k-
 1} \frac{1}{i}$ a.s.

 Indeed, using a Chernoff bound and a large enough $K_0$,
 we obtain that for every $v\in C$
 $$
 \Prob{d (v)>\frac{K_0}{8} pn} \leq  
 e^{-\frac{K_0pn}{8}}\leq n^{-\frac{K_0}{16}}.
 $$
 Therefore we have
 \begin{eqnarray*}
 \Prob{\exists v\in C \sep d(v) > 
 \frac{K_0}{8} pn} \leq
 n\cdot n^{-\frac{K_0}{16}}\rightarrow 0,
 \end{eqnarray*}
 provided $K_0$ is large enough. Then Corollary~\ref{c:C-E} guarantees 
 BoxMaker's win, thus Breaker's win a.s., 
and the proof of Theorem~\ref{t:cg-b} is complete.
 \eopf

Next we give a winning strategy for Maker in the connectivity game,
thus determining the threshold bias $b_{\cal T}^p$ up to a constant factor.

Obviously, Breaker wins if and only if he claims all the edges of
a cut, i.e. all the edges connecting some set of vertices with its
complement. In order to win Maker has to claim one edge in
each of the cuts. This observation enables us to formulate the connectivity
game in a different way, where winning sets are cuts and roles of players
are exchanged -- Breaker wants to occupy a cut and Maker wants to 
prevent Breaker from doing so. To avoid confusion we refer to the
players of this ``cut-game'' by CutMaker and CutBreaker.

This new point of view
enables us to give Maker a winning strategy using Theorem~\ref{t:es_cap}, 
which is a criterion for CutBreaker's win.
Observe, that in this ``cut-game'' CutBreaker (alias Maker) only
cares about occupying the existing edges of a cut, that's why we are
going to look at the family $\cf_p^{\cap}$ instead of $\cf_p$.

\begin{theorem}\label{t:cg-m}
There exists $k_0>0$, so that for $p>\frac{32\log n}{n}$ and
$b\leq k_0 p\frac{n}{\log n}$
Maker wins the random connectivity game 
$(E(K_n)_p, \ct_p, 1, b)$ a.s.
\end{theorem}

\proof
For $b_0=\frac{\log 2}{2}\cdot\frac{n}{\log n}$ 
we are going to prove that the
conditions of Theorem \ref{t:es_cap} are satisfied if
$\cf$ is the set of all cuts in a complete graph 
with $n$ vertices.

On one hand, Beck \cite{be2} showed
$\sum_{k=1}^{n/2} {n \choose k} 2^{-\frac{k(n-k)}{b_0}}\rightarrow 0$,
which means that condition (\ref{c:es}) holds in this setting. 

On the other hand, for a cut $A\in\cf$ we have $|A|\geq n-1$ which implies 
condition (\ref{c:ws_large}). 
If we set $\delta=1/2$ we can apply Theorem \ref{t:es_cap} which gives
that $(\widehat{E(K_n)}_p, \cf_p^\cap, \frac{\log 2}{4}p\frac{n}{\log 
n}, 1)$ is a CutBreaker's win a.s. 
The statement of the theorem  immediately follows.
\eopf

Theorem~\ref{t:cg-b} and Theorem~\ref{t:cg-m} together imply
part $(i)$ of both Theorem~\ref{main1} and \ref{main2}.

\subsection{Hamiltonian cycle game} \label{ss:ham}

Here we investigate the random version of the  
$(1\:b)$ biased game $(E(K_n), \ch, 1,b)$ on 
the complete graph $K_n$, where $\ch$ is the set of all Hamiltonian cycles. 
Maker's goal is to 
occupy all edges of a Hamiltonian cycle, while Breaker wants to
prevent that.
Breaker can obviously win when Maker is not able to
claim a connected graph and thus from Theorem \ref{t:cg-b} we obtain the
following corollary.

\begin{corollary}\label{t:hg-b}
There exists $H_0>0$ so that for every $p\in [0,1]$ 
and $b\geq H_0 p\frac{n}{\log n}$
Breaker wins the random Hamiltonian cycle game $(E(K_n)_p, \ch_p, 1, b)$ a.s.
\end{corollary}

The next Theorem describes Maker's strategy.

\begin{theorem}\label{t:hg-m}
There exists $h_0>0$, so that for $p>\frac{32\log n}{\sqrt{n}}$
and $b\leq h_0 p\frac{\sqrt{n}}{\log n}$
Maker wins the random Hamiltonian cycle game 
$(E(K_n)_p, \ch_p, 1, b)$ a.s.
\end{theorem}

\proof
Maker wins, if at the end of the game the subgraph $G_M$ (containing the
edges claimed by Maker) has connectivity $\kappa(G_M)$ greater or 
equal than independence number $\alpha(G_M)$. Indeed, from
the criterion of Chv\' atal and Erd\H os for Hamiltonicity \cite{ce2}, we
obtain that $G_M$ then contains a Hamiltonian cycle. 

We show that Maker, using only his odd moves, can 
ensure that the connectivity of his graph
at the end of the game is greater then $k=\sqrt{n}/2$ and, 
using his even moves, can make the independence number at 
the end of the game smaller then $k=\sqrt{n}/2$.
In other words we will look at two separate games
where in each of them Maker plays one move against Breaker's $2b$ 
moves. This is a correct strategy, because 
moves of Maker made in one of these games
cannot hurt him in the other.

We first look at the odd Maker's moves.
To ensure that $\kappa(G_M)\geq k$, Maker has to claim one 
edge in every cut of a graph obtained from the initial graph
by removing some $k$ vertices. More precisely, we are going to 
prove the conditions of
Theorem \ref{t:es_cap} for the biased $(b'\:1)$ game, where 
$b'=\frac{\log 2}{2}\cdot\frac{\sqrt{n}}{\log n}$ and
$$
\cf = \bigg\{ \{v_1 v_2 \sep v_1\in V_1,\, v_2\in V_2\} \sep
V(K_n)= V_0 \cupdis V_1 \cupdis V_2,\, |V_0|=k,\, V_1,V_2\not= 
\emptyset \bigg\}.
$$
That is, Maker plays the role of ``CutBreaker'' by trying to
break all the cuts in $\cf$.

Since the size of each of the sets in $\cf$ is at least $n-k-1$ we 
have
$$
\lim_{n\rightarrow\infty} \min_{A\in\cf} \frac{|A|}{b'} = 
\lim_{n\rightarrow\infty}\frac{2\log n (n-\sqrt{n}/2-1)}{\log 2\, \sqrt{n}} =
\infty,
$$
and the 
condition (\ref{c:ws_large}) holds. Next, we have
\begin{eqnarray*}
\sum_{A\in\cf} 2^{-\frac{|A|}{b'}}
&=& \sum_{i=1}^\frac{n-k}{2} {n \choose i} {n-i \choose k} 
2^{-\frac{i(n-i-k)}{b'}} \\
&<& \sum_{i=1}^{k} n^{2k} 2^{-\frac{n-k-1}{b'}}
+ \sum_{i=k+1}^\frac{n-k}{2}
2^{2n-\frac{k(n-2k)}{b'}}\\
&<& k\cdot n^{-\sqrt{n}} +
n\cdot n^{-n}\rightarrow 0,
\end{eqnarray*}
which gives the condition (\ref{c:es}). Therefore, CutBreaker 
(alias Maker) 
wins the game $(\widehat{E(K_n)}_p, \cf^\cap_p, 
\frac{\log 2}{4}p\frac{\sqrt{n}}{\log n}, 1)$ a.s., provided
$p\geq \frac{32 \log n}{\sqrt{n}}$.

In the other part of the game using even moves Maker has to ensure 
that $\alpha(G_M)\leq k=\sqrt{n}/2$. 
That is going to be true if Maker manages to 
claim at least one edge in every clique of $k$ elements. To prove 
that it is possible we again use Theorem \ref{t:es_cap} for a
biased $(b'\:1)$ game with the same value of $b'= \frac{\log 2}{2}\cdot
\frac{\sqrt{n}}{\log n}$. 
But now 
$\cf$ is the family of the edgesets of all cliques of size $k$ and
Maker will play the role of ``CliqueBreaker'' in this game.

We have
$$
\lim_{n\rightarrow\infty} \min_{A\in\cf} \frac{|A|}{b'} = 
\lim_{n\rightarrow\infty}\frac{2\log n {\frac{\sqrt{n}}{2} \choose 
2}}{\log 2\, \sqrt{n}} = \infty,
$$
and the 
condition (\ref{c:ws_large}) is satisfied. It remains to prove that 
the condition (\ref{c:es}) holds.
\begin{eqnarray*}
\sum_{A\in\cf} 2^{-\frac{|A|}{b'}} &=&{n \choose k} 
2^{-\frac{{k\choose 2}}{b'}}<
\left( \frac{ne}{k}2^{- \frac{k-1}{2b'}}\right)^k \\
&<& 2^{-\sqrt{n}} \rightarrow 0. 
\end{eqnarray*}
Therefore, CliqueBreaker 
wins the game $(\widehat{E(K_n)}_p, \cf^\cap_p, 
\frac{\log 2}{4}p\frac{\sqrt{n}}{\log n}, 1)$ a.s., provided
$p\geq \frac{32\log n}{\sqrt{n}}$.

Putting the two parts of the game together we have that Maker 
wins $(E(K_n)_p, \ch_p, 1, \frac{1}{16} p\frac{\sqrt{n}}{\log 
n})$ a.s. \eopf

Combining the statements of Corollary~\ref{t:hg-b} and 
Theorem~\ref{t:hg-m} we obtain part $(iii)$ of both Theorems
\ref{main1} and \ref{main2}.

\subsection{Perfect matching game} \label{ss:pm}

The upper and lower bounds obtained in the previous subsection 
for the threshold bias of the random Hamiltonian cycle game  
are not tight. We firmly believe that our strategy for Maker in that
game is not
optimal. The game we consider next is simpler for Maker, and for that we
are able to obtain bounds optimal up to a constant factor.

Recall that $\cm$ is
the set of all perfect matchings on $K_n$. We will assume that
$n$ is even. 
In the game $(E(K_n), \cm, 1,b)$ Maker's goal is to 
occupy all edges of a perfect matching, while Breaker wants to prevent that.

The following theorem provides the winning strategy in the random
perfect matching game for Maker. 

\begin{theorem}\label{t:mg-m}
There exists $m_0>0$, so that for $p>64 \frac{\log n }{n}$
and $b\leq m_0 p\frac{n}{\log n}$
Maker wins the random perfect matching game 
$(E(K_n)_p, \cm_p, 1, b)$ a.s.
\end{theorem}

\proof 
We can show that Maker can win in a slightly harder game. More
precisely, if the set of vertices of $K_n$ is partitioned into two
sets $A$ and $B$ of equal size before the game starts, we are going to
show that Maker
can claim a perfect matching with edges going only between $A$ and $B$.

For disjoint sets $X,Y\subset V(K_n)$, we define $E(X,Y)$
to be the set of edges between $X$ and $Y$.
Let $\cf$ be a family of sets of edges,
$$
\cf = \{ E(X,Y) \sep \emptyset \not= X \subset A,\, 
\emptyset \not= Y \subset B,\, |X|+|Y|=\frac{n}{2} +1 \}.
$$

Suppose that at the end of the game Maker has not claimed all edges of
any perfect
matching between $A$ and $B$. Hall's 
necessary and sufficient
condition for existence of a perfect matching implies that there
exist sets $X_0 \subset A$ and $Y_0 \subset B$ such that $|X_0|>|Y_0|$
and all edges in $E(K_n)_p \cap 
E(X_0, B\setminus Y_0)$ were claimed by Breaker.

Therefore, in order to win, Maker has to claim at least one edge in
each of the sets from $\cf$, i.e.\ the game $(\widehat{E(K_n)}_p, 
\cf^\cap_p, b, 1)$, which we call {\em Hall},
should be a HallBreaker's win.

To prove that HallBreaker wins we are going to use Theorem
\ref{t:es_cap}. We set $\delta = 1/2$ and 
$b_0=\frac{\log 2}{4}\cdot \frac{n}{\log n} $ .

First we show that condition (\ref{c:es}) holds. We have
\begin{eqnarray*}
\sum_{k=1}^{n/2} {n/2 \choose k} {n/2 \choose n/2 -k+1}
2^{-\frac{k(n/2-k+1)}{b_0}} &<& 2\sum_{k=1}^{\lfloor n/4 \rfloor} 
{n/2 \choose k}^2 2^{-\frac{k(n/2-k+1)}{b_0}} \\
&<& 2\sum_{k=1}^{\lfloor n/4 \rfloor} 
\left( e^{2\log (n/2) - 2\log n  }\right)^k \\
&=& 2\sum_{k=1}^{\lfloor n/4 \rfloor} \left( \frac{1}{4} \right)^k< 1.
\end{eqnarray*}

Since
$$
\lim_{n\rightarrow\infty} \min_{A\in\cf} \frac{|A|}{b_0}>
\lim_{n\rightarrow\infty} \log n = \infty,
$$
the condition (\ref{c:ws_large}) is also satisfied and we can apply
Theorem \ref{t:es_cap} proving that HallBreaker wins the random 
game $(\widehat{E(K_n)}_p, \cf^\cap_p, \frac{\log 2}{8} 
p\frac{n}{\log 
n}, 1)$ a.s., provided $p>64\log n/n$.

This immediately
implies that Maker wins $(E(K_n)_p, \cm_p, 1, b)$ a.s. \eopf

Theorem \ref{t:cg-b} ensures a win for Breaker in the perfect matching
game, if $b> K_0pn/\log n$. 
This, together with the above Theorem~\ref{t:mg-m}
proves part $(ii)$ of Theorems~\ref{main1} and \ref{main2}.

\subsection{Clique game} \label{ss:clique}

Here we look at the random version of the  
$(1\:b)$ biased clique game $(E(K_n), \ckk_k, 1,b)$ on 
a complete graph $K_n$, where $\ckk_k$ is
the set of all cliques of constant size $k$. 
Maker's goal is to 
occupy all edges of a clique of size $k$ while Breaker wants to 
prevent that.

The deterministic clique game was extensively studied by Bednarska and
\L uczak in
\cite{bl1}. They proved a more general result
by determining the order of the threshold bias for the whole
family of games in which Maker's goal is to claim an arbitrary fixed
graph $H$. In this section, we will largely rely on the
constructions and ideas from their paper.

If $\{F_1,\dots ,F_t\}$ is a family of $k$-cliques having two common
vertices, and $e_i\in E(F_i)$, $i=1,\dots ,t$ are distinct edges, 
then we call
the graph $\cup_{i=1}^t F_i$ a {\em $t$-2-cluster} and
the graph $\cup_{i=1}^t (F_i - e_i)$ a {\em $t$-fan}. If furthermore
the $k$-cliques have three vertices in common, then a
$t$-2-cluster is called a {\em $t$-3-cluster} and
a $t$-fan is called a {\em $t$-flower}.
A $t$-fan or a $t$-2-cluster 
is said to be {\em simple}, if the pairwise intersections 
(of any two $k$-cliques) have size exactly $2$.

In order to prevent Maker to occupy a clique $K_k$, Breaker will play two
auxiliary games. In the first one he prevents Maker from occupying a 
3-cluster of constant size.

\begin{lemma} \label{l:3-clusters}
There exists $t=t(k)$, so that for $\eps=\frac{1}{2(k+2)}$,
$p=\omega (n^{-\frac{2}{k+1}})$ and $b>pn^\frac{2(1-\eps)}{k+1}$ 
Breaker wins the game $(E(K_n)_p,
\mbox{$t$-3-clusters}, 1, b)$ a.s.
\end{lemma}

\proof
To apply Theorem~\ref{t:es}, it is enough to check that 
there exists $t$ such that for the random variable
$$
Y:=\sum_{\mbox{\scriptsize $t$-3-cluster $C$ in $G(n,p)$}} (1+b)^{-
e(C)},
$$
$Y<\frac{1}{b+1}$ holds a.s.

We have
$$
\Ex{Y} 
  = \sum_{\mbox{\scriptsize $t$-3-cluster $C$ in $K_n$}}
\left(\frac{p}{1+b}\right)^{e(C)}.
$$
Let $b_1=\frac{b+1}{p}-1$. In~\cite{bl1}, it is shown that there 
exists $t$ for which
$$
\sum_{\mbox{\scriptsize $t$-3-cluster $C$ in $K_n$}}
\left(\frac{p}{1+b}\right)^{e(C)} \leq K_0\frac{1}{b^{1+k_0}_1},
$$
where $k_0,K_0>0$ are constants depending on $k$. This implies 
$\Ex{Y}=o\left( 
\frac{1}{b+1}\right)$, and by Markov inequality we get that
$Y<\frac{1}{b+1}$ a.s. \eopf

During a game, a $t$-fan (or $t$-flower) is
said to be {\em dangerous} if all the $t$ edges missing from the cliques that
make up the $t$-fan are present in the
graph on which the game is played, but not yet
claimed by any of the players.
Note that if at any moment of the game $(E(K_n)_p,
\mbox{$t$-3-clusters}, 1, b)$ Maker claimed a dangerous
$(b+1)t$-flower, then he could win since he
could claim a $t$-3-cluster in his next $t$ moves by simply claiming
missing edges, one by one.
Hence, Lemma \ref{l:3-clusters} implies the following.

\begin{corollary} \label{c:flowers}
There exists $t=t(k)$ so that for 
$\eps=\frac{1}{2(k+2)}$ and $p=\omega(n^{-\frac{2}{k+1}})$, Breaker 
playing
a $(1\: pn^\frac{2(1-\eps)}{k+1})$ game on edges of random graph
$E(K_n)_p$ can make sure that Maker does
not claim a dangerous $\left( pn^\frac{2(1-\eps)}{k+1}t\right)$-flower at any
moment of the game.
\end{corollary}

Next we deal with the second auxiliary game of Breaker; in this game
he prevents the appearance of too many simple $b^\varepsilon$-fans.

\begin{lemma} \label{l:simple_fans}
There exists $C_0>0$, such that for $\eps_1=\frac{1}{6(k+2)}$,
$p\geq n^{-\frac{2}{k+1}}\log^{1/\eps_1} n$, $b>C_0 p
n^\frac{2}{k+1}$ and $s=b^{\eps_1}$
Breaker wins the game $(E(K_n)_p, \mbox{unions of }
\frac{1}{2}{b \choose s}$ $\mbox{simple
$s$-fans},$ $1, b/2)$ a.s.
\end{lemma}

\proof
Let $c_s(n)$ be the number of simple 
$s$-2-clusters contained in $K_n$, and let
$X_s$ be the random variable counting the number of simple
$s$-2-clusters contained in $G(n,p)$. Using the first moment method we
get
$$
\Prob{X_s\geq \Ex{X_s}\log n} \leq \frac{1}{\log n} \longrightarrow 0,
$$
and using this, a.s.\ we have that
\begin{eqnarray*}
& &\sum_{\mbox{\scriptsize dangerous simple} \atop 
\mbox{\scriptsize $s$-fan $C$ in $G(n,p)$}}
(1+b/2)^{-e(C)} \\
&\leq &  
\sum_{\mbox{\scriptsize simple $s$-2-cluster $K$}\atop
\mbox{\scriptsize in $G(n,p)$}} {k\choose 2}^s
(1+b/2)^{-s({k \choose 2}-2)-1} \\
&\leq & {k \choose 2}^s \log n \cdot c_s(n) p^{s({k \choose 2}-1)+1}
2^{sk^2} b^{-s({k \choose 2}-2)-1} \\
&\leq& \log n\cdot C_1^s {n \choose 2}\frac{{n \choose k-2}^s}{s!}
\left(\frac{p}{b}\right)^{s({k \choose 2}-1)+1} b^s\\
&\leq& n^3 \cdot C_1^s n^{(k-2)s}
\left(\frac{1}{C_0n^{\frac{2}{k+1}}}\right)^{s(k+1)(k-2)/2+1} 
\frac{b^s}{s!}\\
&\leq& n^3 \cdot \left(\frac{C_1}{C_0^{{k\choose 2}-1}}\right)^s 
\left(\frac{1}{C_0n^{\frac{2}{k+1}}}\right) 
\frac{b^s}{s!}  < \frac{1}{2}{b\choose s}\frac{1}{b+1},
\end{eqnarray*}
where $C_1=C_1(k)$ is a constant. The last inequality is valid
since $p\geq n^{-\frac{2}{k+1}}\log^{1/\eps_1} n$, and 
for $C_0$ large enough
$\left( C_1/C_0^{{k\choose 2}-1}  \right)^s\leq n^{-5}$.
This enables us to apply Theorem \ref{t:gen-es}, 
and the statement of the lemma is proved.
\eopf

Now we are ready to state and prove the theorem ensuring Breaker's
win in the clique game on the random graph. In the proof, we are 
going to use this result of Bednarska and \L uczak.

\begin{lemma} \cite{bl1} \label{l:ind_sets}
For every $0<\eps <1$ there exists $b_0$ so that every
graph with $b>b_0$ vertices and at most $b^{2-\eps}$ edges has at
least $\frac{1}{2} {b \choose b^{\eps /3}}$ independent sets of size 
$b^{\eps /3}$.
\end{lemma}

\begin{theorem} \label{t:B_clique}
There exists $C_0>0$ so that for $p\geq 
n^{-\frac{2}{k+1}} \log^{6k+12} n$ and
$b\geq C_0 p n^\frac{2}{k+1}$ Breaker wins 
the  random clique game $(E(K_n)_p, (\ckk_k)_p, 1,b)$ a.s.
\end{theorem}

\proof
Breaker will use $b/2$ of his moves to defend ``immediate threats'',
i.e.\ to claim the remaining edge in all $k$-cliques in which Maker 
occupied all but one edge.
In order to be able to do this Breaker must ensure that he never 
has to block more
than $b/2$ immediate threats, that is, there is no dangerous
$b/2$-fan.

He will use his other $b/2$ moves to prevent
Maker from creating a dangerous $(b/2)$-fan. 

From Corollary \ref{c:flowers} 
we get that Breaker can prevent Maker from claiming a dangerous
$f$-flower (where
$f=tpn^\frac{2(1-\eps)}{k+1}$, $\eps=\frac{1}{2(k+2)}$ and $t$ is a
positive constant)
using less than $b/4$ edges per move. On the other hand, from Lemma
\ref{l:simple_fans} we have that if $C_0$ is large enough 
Breaker can prevent Maker from claiming 
$\frac{1}{2}{b/2 \choose s}$ simple $s$-fans using $b/4$ edges per
move, where $s=(b/2)^{\eps/3}$.

Suppose that Maker managed to claim a dangerous 
$(b/2)$-fan. We define an auxiliary graph
$G'$ with the vertex set 
being the set of all $b/2$ $k$-cliques of this dangerous fan, and
two $k$-cliques being connected with an edge if they have at least $3$
vertices in common. Since there is no dangerous 
$f$-flower in Maker's graph, the
degree of each of the vertices of the graph $G'$ is at most $fk$ and
therefore $e(G')<\frac{bfk}{2}\leq
\left(\frac{b}{2}\right)^{2-\eps}$. 
On the other hand, the number of independent sets
in $G'$ of size $s$ cannot be more than $\frac{1}{2}{b/2 \choose
s}$, since each of the independent sets in $G'$ corresponds to 
a simple $s$-fan in Maker's graph.

Since the last two facts are obviously in contradiction with Lemma
\ref{l:ind_sets}, 
Maker cannot claim a dangerous $b/2$-fan and the statement of the
theorem is proved. \eopf

To prove the theorem for Maker's win, we need the
following lemma which is a slight modification of a result from 
\cite{bl1}. Let
$G(n,M)$ denote
the graph obtained by choosing a graph on $n$ vertices with $M$ edges
uniformly at random.

\begin{lemma} \label{l:bl1}
There exists $0<\delta_k<1$, such that for $M=2
\lfloor n^{2-2/(k+1)} \rfloor$ a.s.\ each
subgraph of $G(n,M)$ with $\lfloor (1-\delta_k) M\rfloor$ edges contains
a copy of $K_k$.
\end{lemma}

\proof
For $0<\delta_k<1$, we call a subgraph $F$ of 
$K_n$ bad, if $F$ has $M$ edges and it contains a subgraph $F'$ 
with 
$\lfloor (1-\delta_k) M\rfloor$ edges that does not contain a copy of 
$K_k$.
In \cite{bl1}, it is proved that there exist constants 
$0<\delta_k<1$ and $c'_1>0$
such that the number of bad subgraphs of 
$K_n$ is bounded from above by
$$
e^{-c'_1 M/6} {{n \choose 2}\choose M } = o(1) {{n \choose 
2}\choose M }.
$$
\eopf

Using the last lemma 
we can prove a theorem for Maker's win in the random clique game.

\begin{theorem} \label{t:M_clique}
There exists $c_0>0$ so that for $p>\frac{1}{c_0}
n^{-\frac{2}{k+1}}$ and
$b\leq c_0 p n^\frac{2}{k+1}$ Maker wins 
the  random clique game $(E(K_n)_p, (\ckk_k)_p, 1,b)$ a.s.
\end{theorem}

\proof
We will follow the analysis of the random Maker's strategy 
proposed in \cite{bl1}, looking at
$G(n,M')$, where $M'=p{n\choose 2}$.
We will prove that the $k$-clique game on $G(n,M')$ is a Maker's
win a.s., which implies that 
the same is true on $G(n,p)$, as being a Maker's
win is a monotone property \cite[Chapter 2]{bol}. 

In each of his moves Maker chooses one of the edges
of $G(n,M')$ that was not previously claimed by him, uniformly at
random. If the edge is free he claims it and we call that a
successful Maker's move. If the 
edge was already claimed by Breaker, then Maker skips
his move (e.g.\ claims
an arbitrary free edge, and that edge we will not encounter for the
future analysis).

Let $0<\delta_k<1$ be chosen so that the conditions
of Lemma \ref{l:bl1} are satisfied.
We look at the course of game after $M=2\lfloor n^{2-2/(k+1)} \rfloor$ moves.

By choosing $c_0\leq \delta_k/12$, we have
\begin{eqnarray*}
M &\leq& \frac{\delta_k}{6 c_0} \lfloor n^{2-2/(k+1)} \rfloor \\
&\leq& \frac{\delta_k}{2} \frac{1}{b+1} p{n \choose 2}.
\end{eqnarray*}
That means that only at most $\delta_k /2$ fraction of the total 
number of
elements of the board $E(G(n,M'))$ is claimed (by both players)
after move $M$. Therefore,
the probability that the edge randomly chosen in Maker's $m$th 
move,
$m\leq M$, is already claimed by Breaker is bounded from above by
$\delta_k/2$. That means that Maker 
has at least $(1-\delta_k)M$ successful moves a.s. 

Since in each of his moves 
Maker has chosen edges uniformly at random (without repetition) 
from $E(G(n,M'))$, the graph
containing edges chosen by Maker in his first $M$ moves (both
successful and unsuccessful)
actually is a random graph $G(n,M)$. 
Applying Lemma \ref{l:bl1}, we get that the graph
containing edges claimed by Maker in his successful moves contains a
clique of size $k$ a.s., which means that a.s.\ there
exists a non-randomized winning strategy for Maker.
\eopf

Combining the statements of Theorem~\ref{t:B_clique} and
Theorem~\ref{t:M_clique} we obtain part $(iv)$ of Theorem~\ref{main2}.

\section{Unbiased games} \label{s:unbiased}

\subsection{Connectivity one-on-one} \label{ss:con1on1}

A theorem of Lehman enables us to determine the threshold probability
$p_{\ct}$ with extraordinary precision. Namely,
Lehman \cite{L} proved that the unbiased connectivity game is won by
Maker (now as a second player!) if and only if 
the underlying graph contains two edge-disjoint spanning trees.
The threshold for the appearance of two edge-disjoint spanning trees was
determined exactly by Palmer and Spencer~\cite{PS}. 

To formulate the consequence of these two results we need the concept of
{\em graph process}. Let $e_1,\ldots e_m$ be the edges of $K_n$, where
$m={n\choose 2}$. Choose a permutation $\pi\in S_m$ uniformly at random
and define an increasing sequence of subgraphs $(G_i)$ where $V(G_i)=V(K_n)$
and $E(G_i)=\{e_{\pi(1)},\ldots , e_{\pi (i)}\}$. It is clear that
$G_i$ is an $n$-vertex graph with $i$ edges, selected uniformly at random
from all $n$-vertex graphs with $i$ edges.

Given a particular graph process $(G_i)$ and a graph property
$\cp$ possessed by $K_n$, the {\em hitting time} 
$\tau (\cp)=\tau (\cp, (G_i))$ is the minimal $i$ for which 
$G_i$ has property $\cp$.

The consequence of the theorems of Lehman, and Palmer and Spencer is
that the very moment the last vertex receives its second adjacent edge,
the unbiased connectivity game is won by Maker a.s.
More precisely, the following is true.
\begin{corollary} \label{connectivity-precise}
For the unbiased connectivity game we have that a.s.
$$\tau(\mbox{Maker wins $\ct$})=\tau(\mbox{$\exists$ two edge-disjoint spanning trees})=
\tau (\delta(G)\geq 2).$$ 
\end{corollary}

In particular,
for edge-probability $p=(\log n + \log\log n +g(n))/n$, where 
$g(n)$ tends to infinity arbitrarily slowly, Maker wins the unbiased
connectivity game a.s., while if $g (n)\rightarrow -\infty$, then
Breaker wins a.s. 

\heading{Remark.} The assumption that Maker is the second player is just
technical, for the sake of smooth applicability of Lehman's Theorem.
If Maker is the first player, then from the proof of
Lehman's Theorem one can infer that Maker wins if and only if the base graph
contains a spanning tree and a spanning forest of two components, which are
edge-disjoint. This property has the same sharp threshold as 
the presence of two edge-disjoint spanning trees, and the hitting time
should be the same when the next to last vertex receives its second
incident edge.

\subsection{$k$-cliques one-on-one} \label{ss:clique1on1}

Let us fix $k$ and
let $(F_1,\ldots , F_s)$ be a sequence of $k$-cliques. Then 
$F=\cup_{i=1}^sF_i$ is called an {\em $s$-bunch} if 
$V(F_i)\setminus (\cup_{j=1}^{i-1} V(F_j))\not= \emptyset$ and
$|V(F_i)\cap (\cup_{j<i}V(F_j))|\geq 2$, for each $i=2,\ldots  , s$.
Recall that 
an $s$-bunch in which the pairwise intersection of any two cliques 
is the same two vertices, was called a {\em simple $s$-$2$-cluster}. Let us denote 
the simple $s$-$2$-cluster by $C_s$. 

For a graph $G$, the {\em density} of $G$ is defined as
$d(G)=\frac{e(G)}{v(G)}$, and
the {\em maximum density} of $G$ is defined as
$m(G)=\max_{H\subseteq G} d(H)$. A graph $G$ with $m(G)=d(G)$ is called
{\em balanced}. The maximum density of a graph $G$ determines the threshold
probability for the appearance 
of $G$ in the random graph.
More precisely, $(i)$ if $p=o(n^{-1/m(G)})$, then $G(n,p)$ does not
contain $G$ a.s., and $(ii)$ if 
$p=\omega(n^{-1/m(G)})$, then $G(n,p)$ does contain $G$ a.s.

We need two properties of simple $s$-$2$-clusters and $s$-bunches.

\begin{lemma} \label{l:mcs}
For every positive integer $s$, 
$C_s$ is balanced and has maximum density
$m(C_s)=d(C_s)=\frac{k+1}{2}-\frac{k}{sk-2s+2}$. 
\end{lemma}

\proof
It is easy to check that $v(C_s)=s(k-2)+2$,
$e(C_s)=s{k\choose 2} -s+1$, and thus
$d(C_s)=\frac{e(C_s)}{v(C_s)}=\frac{k+1}{2}-\frac{k}{sk-2s+2}$.

Let $T$ be a subgraph of $C_s$. We want to prove $d(T)\leq d(C_s)$.
Since $C_s$ is the union of $k$-cliques, 
$C_s=\cup_{i=1}^s F_i$, if we set $E_i=F_i\cap T$ 
we have that $T=\cup_{i=1}^s E_i$, and we can assume that
each $E_i$ is a clique of order $k_i\leq k$. 

We can also assume that the two vertices in $\cap_{i=1}^s V(F_i)$ are
in $T$, since otherwise their inclusion would 
increase the density. This implies $k_i\geq 2$ for $i=1,\ldots ,s$. 

Let us relabel the cliques in such a way that $E_i\neq F_i$ if and only if
$i=1,\ldots , s_1$. Then
$$\frac{e(C_s)}{v(C_s)} \geq \frac{e(T)}{v(T)}=
\frac{e(C_s)-\sum_{i=1}^{s_1} \left( {k\choose 2}
-{k_i\choose 2}\right)}{v(C_s)-\sum_{i=1}^{s_1} (k-k_i)},
$$
since
$$\frac{e(C_s)}{v(C_s)}<\frac{k+1}{2}\leq
\frac{\sum_{i=1}^{s_1} \left( k-k_i\right)
\frac{k+k_i-1}{2}}{\sum_{i=1}^{s_1} (k-k_i)}.
$$
The last inequality is true since the last fraction is
the weighted average
of the numbers $(k+k_i-1)/2$, each of them being at least $(k+1)/2$.
\eopf

\begin{lemma}\label{megegy}
Let $s\geq 3$ be a positive integer. No $s$-bunch has smaller maximum 
density than the simple $s$-$2$-cluster.
\end{lemma}

\proof
When $k=3$, the $s$ bunch is a union of triangles. 
Then any $s$-bunch has the same number of vertices as the simple 
$s$-$2$-cluster, while the number of edges, and thus the density is minimized for
the simple $s$-$2$-cluster.

From now on let us assume that $k\geq 4$.
Let $s\geq 3$, and let $(F_1, F_2,\dots , F_s)$ be the
sequence of $k$-cliques of an arbitrary $s$-bunch $B_s=\cup_{i=1}^s F_i$. 
For every $i\in\{2,3,\dots ,s\}$, let $F'_i= \left( \cup_{j=1}^{i-1}
F_j \right) \cap F_i$.
Then, we have
\begin{eqnarray*}
d(B_{s}) &=& \frac{s{k \choose 2} - \sum_{i=2}^s e(F'_i)}{ 
sk-\sum_{i=2}^s v(F'_i)} \\
&=& \frac{e(C_s) - \sum_{i=2}^s (e(F'_i) -1)}{ 
v(C_s)-\sum_{i=2}^s (v(F'_i) -2)} \\
&\geq & \frac{e(C_s) - \sum_{i=2}^s ({v(F'_i) \choose 2} -1)}{ 
v(C_s)-\sum_{i=2}^s (v(F'_i) -2)} \\
&\geq & \frac{e(C_s)}{v(C_s)}.
\end{eqnarray*}

In the last inequality the terms with $v(F_i')=2$ disappear, and
otherwise we use that $v(F'_i)\leq k-1$ for every $i$, so
$\frac{{v(F'_i) \choose 2} -1}{v(F'_i) -2} \leq \frac{k}{2} \leq 
\frac{e(C_s)}{v(C_s)}$.

Hence, simple $s$-$2$-clusters 
have the smallest density among all $s$-bunches.
For any $s$-bunch $B_s$ and the simple $s$-$2$-cluster 
$C_s$ we immediately obtain 
$$
m(B_s)\geq d(B_s)\geq d(C_s) = m(C_s),
$$
and the lemma is proved.
\eopf
\heading{Remark.} The previous lemma is of course true for $s=1$, but not for
$s=2$.

As a consequence of the last two lemmas we get a strategy for
Breaker in the $(1:1)$ clique game.

Let $H$ be a graph and consider the auxiliary graph
$G_H$ with vertices corresponding to the $k$-cliques of $H$,
two vertices being adjacent if the corresponding cliques 
have at least two vertices in common. 
Let $F_1,\ldots, F_s$ be the cliques corresponding to a connected
component of $G_H$. Then the graph $\cup_{i=1}^sF_i$ is called 
an {\em $s$-collection} or just a {\em collection} of $H$. 
Note that the edgeset of any $H$ is uniquely partitioned into sets
$N$ and $E(A_i)$, where $N$ contains the edges which do not participate in
a $k$-clique, while the $A_i$ are the collections of $H$.

\begin{theorem} \label{p:kk_1on1}
For every $k\geq 4$ and $\varepsilon >0$, 
$p_{\ckk_k} \geq n^{-\frac{2}{k+1}-\varepsilon}$.
For $k=3$, we have that
$p_{\ckk_3} \geq n^{-\frac{5}{9}}$.
\end{theorem}

\proof
First we give a strategy for Breaker to win $\ckk_k$ if the game is played
on the edgeset of a $(2k-4)$-degenerate graph $L$. 
Consider the ordering $v_1,\ldots , v_{v(L)}$ of $V(L)$,
such that $|N_{V_j}(v_{j+1})|\leq 2k-4$ for $j=1,\ldots , v(L)-1$, where 
$V_j=\{ v_1,\ldots , v_j\}$. Then Breaker's strategy is the following:
if Maker takes an edge connecting $v_{j+1}$ to $V_j$, 
then Breaker takes another one also connecting $v_{j+1}$ to $V_j$. 
If there is no such edge available, then
Breaker takes an arbitrary edge. Suppose for a contradiction that
Maker managed to occupy a $k$-clique $v_{i_1}, \ldots, v_{i_k}$ 
against this strategy, where $i_1 < \cdots < i_k$.
This is impossible, since Maker could have never 
claimed $k-1$ of the edges $v_{j}v_{i_k}$, $j<i_k$.

Let $E(K_n)_p=N\cupdis E(A_1)\cupdis\dots \cupdis E(A_{h})$ be the
partition of the edges, 
such that $N$ contains all edges that do not participate in
any $k$-clique, and each $A_i$ is a collection of $k$-cliques. 
(Corresponding to the connected
components of the auxiliary graph $G_{G(n,p)}$
defined on the set of $k$-cliques of $G(n,p)$.)

Breaker can play the game $(E(K_n)_p, (\ckk_k)_p,1,1)$ by playing
separately on each of the sets $E(A_i)$. More precisely, whenever
Maker claims an edge which is in some $E(A_i)$, Breaker can play according
to a strategy restricted just to $E(A_i)$. Since, crucially, 
the edgeset of
each $k$-clique is completely contained in exactly one of the $E(A_i)$, 
Maker can only
win the game on $E(K_n)_p$ if he wins on one of the $E(A_i)$.

Now we are going to show that every collection $A$ on
$v(A)=v$ vertices contains a $\lceil \frac{v-2}{k-2} \rceil$-bunch.
We take an arbitrary $k$-clique $F_1$
from $A$, and build a bunch recursively as follows. If we picked
$k$-cliques $F_1, \dots, F_i$ then we choose $F_{i+1}$ such that
$|V(F_{i+1})\cap (\cup_{j=1}^i V(F_j) )| \geq 2$ and 
$V(F_{i+1})\setminus (\cup_{j=1}^i V(F_j) ) \not= \emptyset$. Note
that this means that $\cup_{j=1}^{i+1} F_j$ is an 
$(i+1)$-bunch. Since the auxiliary graph $G_{A}$ of the collection 
is connected we can keep doing this until $V(A)=\cup_{j=1}^{i_0}
V(F_j)$ for some $i_0$.
Knowing that $v(F_i)=k$ for all $i\leq i_0$, we have
$i_0\geq 1+ \frac{v-k}{k-2} = \frac{v-2}{k-2}$.
So there exists an $\lceil\frac{v-2}{k-2}\rceil$-bunch which is a
subgraph of $A$.

We first look at the case $k\geq 4$.
Let $\varepsilon>0$ be a constant. From 
Lemma~\ref{l:mcs} it follows that there exists an integer $v$ 
such that for $s_0=\lceil\frac{v-2}{k-2}\rceil$ we have
$m(C_{s_0})\geq \frac{k+1}{2}-\frac{k}{v} 
> \left( \frac{2}{k+1} + \varepsilon \right)^{-1}$.
Then for $p = O( n^{-\frac{2}{k+1}-\varepsilon})$
it follows that there is no $s_0$-bunch in $G(n,p)$ a.s.,
since we have have that the first $s_0$-bunch 
that appears in the random graph is the one of the minimum maximum
density, which, by Lemma~\ref{megegy}, is the simple 
$s_0$-$2$-cluster.
Note here that there is a constant (depending on $k$ and $\varepsilon$) 
number of nonisomorphic $s_0$-bunches.

Since in $G(n,p)$ there are no $s_0$-bunches a.s., 
there are also no collections on $v$ vertices a.s. 

Finally, all the collections $A_i$ are
$(2k-4)$-degenerate a.s., since graphs which are not
$(2k-4)$-degenerate have maximum density at least
$\frac{2k-3}{2}\geq \frac{k+1}{2}$, provided $k\geq 4$. 
Note that we know already that a.s. all collections have order at most 
$v$ and thus there are at most a constant (depending on $k$ and $\epsilon$)
number of nonisomorphic non-$(2k-4)$-degenerate graphs.

This proves that Breaker has a winning strategy a.s., if
$k\geq 4$ and $p = O( n^{-\frac{2}{k+1}-\varepsilon})$.

Next, we look at the case $k=3$. As we saw,
any collection of triangles on $v$ vertices contains a
$(v-2)$-bunch. Thus
for $p=o(n^{-5/9})$, no $v$-collection with $v\geq 15$ will appear in
$G(n,p)$ a.s., since it would contain a $13$-bunch, whose
maximum density is at least $m(C_{13})= 2 -\frac{3}{15}$.
This observation makes the problem finite: one has to check
who wins on collections up to $14$ vertices.

Suppose that Maker can win the triangle game on some collection of
triangles on $v\leq 14$ vertices and with maximum density less then
$9/5$. Let $A$ be a minimal such collection (Maker cannot win on
any proper subcollection of $A$). 

If there was
a vertex $w\in V(A)$ with $d_A(w)\leq 2$, the minimality of $A$ would
imply that Breaker has a winning strategy on $A$. Indeed, Breaker plays
according to his strategy on $A-w$, 
and as soon as Maker claims one edge adjacent
to $w$ Breaker claims the other edge adjacent to $w$ (if that edge exists 
otherwise he does not move). This would mean that
Breaker can win on $A$, a contradiction.
Thus, $\delta_A \geq 3$.

Let $B$ be a $(v-2)$-bunch contained in $A$, with $V(A)=V(B)$.
Since $\delta_B = 2$, we have $e(A)\not=
e(B)$. Then
$$
2-\frac{3}{v}=m(C_{v-2})\leq \frac{e(B)}{v} < \frac{e(A)}{v}
<\frac{9}{5},
$$
and 
$$
2v-3 = e(C_{v-2}) \leq e(B) < e(A) < \frac{9v}{5}.
$$
It is easy to check that Maker  cannot win the
game on a graph with less then $5$ vertices, thus $v>4$,
so $e(B)=e(C_{v-2})$ and $e(A)-e(B) = 1$.

Let $\{e\} = E(A)\setminus E(B)$, and let $T_1,\dots, T_{v-2}$ be
the sequence of triangles whose union is the $(v-2)$-bunch $B$. Since
$e(B)=e(C_{v-2})$, for every $i=2, \dots, v-2$ we have that $T_i$ has
a common edge with $\cup_{j=1}^{i-1} T_j$. Then 
$B$ must have at least $2$
vertices of degree $2$. From $\delta_{B\cup\{e\}} = \delta_{A} =3$
we obtain that $B$ has exactly two vertices $b_1,b_2$ with
$d_B(b_1)=d_B(b_2)=2$, and moreover $e=\{b_1, b_2\}$.
Since $e$ has to participate in at least one 
triangle of the collection
$A$, $b_1$ and $b_2$ have to be connected with a $2$-path in $B$,
which is possible only if 
all $T_1,\dots,T_{v-2}$ share a vertex. That means that
$A$ is a $(v-1)$-wheel and it is easy to see that Breaker can win
the triangle game on a wheel of arbitrary size by a simple pairing
strategy. 

This contradiction proves that
for $p=o(n^{-5/9})$, a.s.\ there is no triangle 
collection in $G(n, p)$ on which Maker can win, 
which means that Breaker a.s.\ wins the game on the whole graph.
\eopf

From Theorem~\ref{t:M_clique} we get that Maker can
win the game $(E(K_n)_p, (\ckk_k)_p,1,1)$ for $p=\Theta
(n^{-\frac{2}{k+1}})$ and thus we immediately obtain $p_{\ckk_k}=O
(n^{-\frac{2}{k+1}})$. For the triangle game $\ckk_3$ a 
stronger upper bound can be found.

\begin{proposition} \label{p:triangle}
The game $(E(K_n)_p, (\ckk_3)_p,1,1)$ is a Maker's win a.s.,
provided $p=\omega(n^{-\frac{5}{9}}$).
\end{proposition}

\proof
It is easy to check that Maker can claim a triangle in 
the $(1\: 1)$ game if the board on which the game is played is 
$K_5$ minus an edge. Therefore, as soon as the
graph $G(n,p)$ contains $K_5-e$ a.s., the initial 
game can be won by Maker a.s. 
\eopf

Theorem~\ref{t:M_clique}, Theorem~\ref{p:kk_1on1} 
and Proposition \ref{p:triangle} 
imply parts $(iv)$ and $(v)$ of Theorem~\ref{main1}.

\section{Open questions} \label{s:open}

\heading{More sharp thresholds?} We saw in the previous section 
that the connectivity game has a sharp threshold, and even
more. We think that both the perfect matching game and the 
Hamiltonian cycle game have the same sharp threshold $\frac{\log n}{n}$, 
and maybe even more\dots 
It would be very interesting to decide whether the following 
conjectures are true.

\begin{conjecture}
\begin{eqnarray*}
& (i) & \tau(\mbox{Maker wins $\cm$}) = \tau(\delta(G)\geq 2),\\
& (ii) & \tau(\mbox{Maker wins $\ch$}) = \tau(\delta(G)\geq 4).
\end{eqnarray*}
\end{conjecture}

\heading{Clique game/$H$-game.} 
The exact determination of the threshold $p_{\ckk_k}$ for the
$k$-clique game remains outstanding.

\begin{problem}
Decide whether $p_{\ckk_k}=n^{-\frac{2}{k+1}}$ for $k\geq 4$.
\end{problem}

The arguments of Bednarska and \L uczak  \cite{bl1} could be
extended to full generality to positional games on random graphs
along the lines of Section~\ref{ss:clique}. 
More precisely, the following is 
true. Let $\ckk_H$ be the family of subgraphs of $K_n$, isomorphic to $H$.
Then for any fixed graph $H$ there is a constant $c(H)$, such that
$$b_{\ckk_H}^p=\Theta\left (pb_{\ckk_H}\right) = \Theta \left(p n^{-1/m'(H)}
\right),$$ provided $p\geq \Omega\left( \frac{\log ^{c(H)} n}{n^{1/m'(H)}}
\right).$ 

Concerning the one-on-one game, it would be desirable to determine those 
graphs for which an extension of the low-density Maker's win,
\`{a} la Proposition~\ref{p:triangle}, exists.

\begin{problem}
Characterize those graphs $H$ for which there exists a constant 
$\epsilon(H) >0$, such that
the unbiased game $\ckk_H$ is a.s.\ a Maker's win 
if $p=n^{-1/m'(H)-\epsilon (H)}$.
\end{problem}

For such graphs the determination of the
threshold $p_{\ckk_H}$ is a finite problem, in a way similar to the
case $H=K_3$.

\heading{Relationships between thresholds.} 
It is an intriguing task to understand under what
circumstances the following is true.

\begin{problem} \label{p:gen1}
Characterize those games $(X, \cf)$ for which 
$$p_\cf=\frac{1}{b_\cf}.$$
More generally, characterize the games for which
$$
b_\cf^p = \Theta \left(pb_\cf\right),
$$
for every $p = \omega \left( \frac{1}{b_\cf}\right)$. 
\end{problem}

This is not true in general as the triangle game shows. What is the 
reason it is true for the connectivity game and the perfect matching
game? 
Is it because the appearance of these properties has a sharp
threshold in $G(n,p)$? Or because the winning sets are not of constant size?

\begin{problem} \label{p:gen2}
Suppose $p_\cf=1/b_\cf$. Is it true that for every $p\geq p_\cf$,
$b_{\cf}^p=\Theta \left( pb_\cf\right)$?
\end{problem}

It would be very interesting to relate the thresholds $b_\cf$ and $p_\cf$
to some thresholds of the family $\cf$ in the random graph $G(n,p)$ (or, more
generally, in the random set $X_p$).
It seems to us that if the family $\cf_p$ is quite dense and well-distributed
in $X$, then Maker still wins the $(1\:1)$ game. 

\begin{problem}
Characterize those games $(X,\cf )$ for which there exists a constant 
$K$, such that for any probability $p$
with $\Prob{\min_{x\in X_p}|\{ F\in \cf_p : x\in
  F\}|>K}\longrightarrow 1$,
we have $p_\cf =O(p)$ and/or
$b_\cf =\Omega (1/p)$.
\end{problem}

\paragraph{Acknowledgments.}
We would like to thank the anonymous referees for their thorough
work. Their numerous suggestions and simplifications
greatly improved the presentation of the paper. We are also 
indebted to Ma\l gorzata
Bednarska and Oleg Pikhurko whose remarks about 
Theorem~\ref{p:kk_1on1} not only shortened the proof, but also improved 
the statement.

\end{document}